\documentclass[10pt,a4paper]{article}
\usepackage{graphicx}
\usepackage{amsmath,amsfonts,amsthm}
\usepackage{xfrac,esint,mathrsfs,color}
\usepackage[utf8]{inputenc}
\usepackage[colorlinks=true]{hyperref}
\usepackage[margin=3cm]{geometry}

\newcommand{\Huno}{{\calH}^{1}}
\newcommand{\Rdue}{\R^2}

\newcommand{\C}{\mathbb{C}}
\newcommand\res{\mathop{\hbox{\vrule height 7pt width .5pt depth 0pt
\vrule height .5pt width 6pt depth 0pt}}\nolimits}
\newcommand{\Ln}{{\calL}^2}

\newcommand\R{\mathbb{R}}

\newcommand\calH{\mathcal{H}}
\newcommand\calF{\mathcal{F}}

\newcommand\calA{\mathcal{A}}

\newcommand\calL{\mathcal{L}}
\newcommand{\LM}[1]{\hbox{\vrule width.2pt \vbox to#1pt{\vfill \hrule width#1pt height.2pt}}}
\newcommand{\LL}{{\mathchoice{\,\LM7\,}{\,\LM7\,}{\,\LM5\,}{\,\LM{3.35}\,}}}

\newcommand{\be}[1]{\begin{equation} #1 \end{equation}}

\newtheorem{theorem}{Theorem}[section]

\newtheorem{proposition}[theorem]{Proposition}

\definecolor{green}{rgb}{0,.5,0}

\numberwithin{equation}{section}
\newcounter{Nummer}

\usepackage{fancyhdr}
\begin{document}
\begin{center}
  {\Large
Existence of minimizers for the $2$d stationary Griffith fracture model}\\[5mm]
{\today}\\[5mm]
Sergio Conti$^{1}$, Matteo Focardi$^{2}$, and Flaviana Iurlano$^{1}$\\[2mm]
{\em $^{1}$
 Institut f\"ur Angewandte Mathematik,
Universit\"at Bonn\\ 53115 Bonn, Germany}\\[1mm]
{\em $^{2}$ DiMaI, Universit\`a di Firenze\\ 50134 Firenze, Italy}\\[3mm]
    \begin{minipage}[c]{0.8\textwidth}
    We consider the variational formulation of the Griffith fracture model in two spatial dimensions and prove
existence of strong minimizers, that is deformation fields
which are continuously differentiable
outside a closed jump set and which minimize the relevant energy. To this aim, we show that minimizers of the weak formulation 
of the problem, set in the function space $SBD^2$ and for which existence is well-known, are actually strong 
minimizers following the approach developed by De Giorgi, Carriero, and Leaci in the corresponding scalar setting of the 
Mumford-Shah problem.
    \end{minipage}
\end{center}

\section{Introduction}
The study of brittle fracture in solids is based on the Griffith model, which combines
elasticity with a term proportional to the surface area opened by the fracture. In its variational formulation 
the energy
\begin{equation}\label{eqgriffintro}
 E[u,\Gamma]:=\int_{\Omega\setminus\Gamma} \Big(\frac12 \C e(u)\cdot e(u)+g(u)\Big) dx + \gamma \calH^{n-1}(\Gamma) 
\end{equation}
is minimized over all closed 
sets $\Gamma\subset\Omega$ and all deformations $u\in C^1(\Omega\setminus\Gamma,\R^n)$ subject to suitable
boundary and irreversibility conditions. Here $\Omega\subset\R^n$ is the reference configuration, the function 
$g\in C^0(\R^n)$ represents external volume forces, $e(u)=(\nabla u+\nabla u^T)/2$ is the elastic strain,
$\C\in \R^{(n\times n)\times (n\times n)}$ is the matrix of elastic coefficients, $\gamma>0$ the surface energy.
The evolutionary problem of fracture can be modeled as a sequence of variational problems, in which
one minimizes (\ref{eqgriffintro}) subject to varying loads with a kinematic restriction representing
the irreversibility of fracture, see \cite{FrancfortMarigo1998,BourdinFrancfortMarigo2008} and references therein.

Mathematically, (\ref{eqgriffintro}) is a vectorial free discontinuity problem. Much better understood is its scalar 
version, in which one replaces the elastic energy by the Dirichlet integral,
\begin{equation}\label{eqgMSintro}
 E_\mathrm{MS}[u,\Gamma]:=\int_{\Omega\setminus\Gamma}\Big(\frac12 |Du|^2+g(u)\Big) dx+\gamma\calH^{n-1}(\Gamma) \,,
\end{equation}
and one minimizes over all closed sets $\Gamma\subset\Omega$ and maps $u\in C^1(\Omega\setminus\Gamma,\R)$.
This scalar reduction coincides with the Mumford-Shah functional of image segmentation, that has been widely studied analytically and 
numerically \cite{mum-shah,amb-tort1,AmbrosioFuscoPallara,David2005}. 
By taking into account the structure of the energy (\ref{eqgMSintro})  it is natural to introduce the space $SBV(\Omega)$ of special functions of bounded 
variation, by imposing that the distributional derivative $Du$ is a bounded measure, i.e. $u\in BV(\Omega)$, that can be written as 
$Du=\nabla u \calL^n\res\Omega + [u]\nu_u \calH^{n-1}\LL J_u$ with $\nabla u$ the approximate gradient of $u$, $[u]$ the jump of 
$u$, $J_u$ the $(n-1)$-rectifiable jump set of $u$, $\nu_u$ its normal. 
Therefore, the relaxation of (\ref{eqgMSintro}) is 
\begin{equation}\label{e:weakMS}
 E_\mathrm{MS}^*[u]:=\int_{\Omega} \Big(\frac12 |\nabla u|^2 dx+g(u)\Big) + \gamma \calH^{n-1}(J_u) \,,
\end{equation}
and it is finite provided $u$ belongs to the subspace $SBV^2(\Omega)$ of functions in $SBV(\Omega)$ with
approximate gradient $\nabla u\in L^2(\Omega;\R^n)$ and $\calH^{n-1}(J_u)<\infty$. 
Existence of minimizers for the relaxed problem $E_\mathrm{MS}^*$ follows then from the general compactness 
properties of $SBV^2$, see \cite{AmbrosioFuscoPallara} and references therein.

The breakthrough in the quest for an existence theory for
 the Mumford-Shah functional (\ref{eqgMSintro})
came with the proof by De Giorgi, Carriero and Leaci in 1989 \cite{DegiorgiCarrieroLeaci1989}
that the jump set of minimizers $u$ is essentially closed, in the sense that
 minimizers of the relaxed functional $ E_\mathrm{MS}^*$ obey 
\begin{equation}\label{eqjuclosedintro}
\calH^{n-1}(\overline {J_u}\setminus J_u)=0,\quad\text{or equivalently } \calH^{n-1}(J_u)=\calH^{n-1}(\overline {J_u}).  
\end{equation}
From this, elliptic regularity implies then that $(u,\overline{J_u})$ is a minimizer of the functional in the original formulation (\ref{eqgMSintro}).

We address here the analogous existence issue for (\ref{eqgriffintro}) 
in two spatial dimensions. We assume that $\C$ is a symmetric linear map from $\R^{n\times n}$ to itself 
with the properties
\begin{equation}\label{eqassC}
 \C(z-z^T)=0 \text{ and } \C z\cdot z \ge \alpha |z+z^T|^2 \text{ for all } z\in\R^{n\times n}
\end{equation}
for some $\alpha>0$.
 Our main result is the following
\begin{theorem}\label{theop2}
 Let $\Omega\subset\R^2$ be a bounded Lipschitz set, $g\in C^0(\R^2)$, $\C$ obey the positivity 
 condition (\ref{eqassC}), {$M>0$}. Then the functional (\ref{eqgriffintro}) has a minimizer in the class
 \begin{equation}
  \calA:=\{(u,\Gamma): \Gamma\subset\overline \Omega \text{ closed, } u\in C^1(\Omega\setminus\Gamma,\R^2)\cap L^\infty(\Omega,\R^2), \ {\|u\|_{L^\infty(\Omega,\R^2)}\leq M}\}.
 \end{equation}
\end{theorem}
{The $L^\infty$-bound on deformations  guarantees existence of solutions to the corresponding relaxed problem.}

The proof is sketched below and will be discussed in detail elsewhere \cite{ContiFocardiIurlano2016b}. In \cite{ContiFocardiIurlano2016b} we 
also consider generalizations of the basic model (\ref{eqgriffintro}) with 
$p$-growth, $p\in(1,\infty)$, which may be appropriate for the study of materials with defects, such as damage 
or dislocations,  and are obtained by replacing the quadratic volume energy density with 
\begin{equation}\label{eqdeffintro}
 f_\mu(\xi):=\frac 1p\left(\big(\mathbb{C}\xi\cdot\xi+\mu\big)^{\sfrac p2}-\mu^{\sfrac p2}\right)
\end{equation}
where $\mu\geq 0$ is a parameter. 

This result is restricted to the two-dimensional case because the approximation result in 
Proposition \ref{prop:ricopr} below is only valid in two dimensions.

\section{Outline of the proof}

Following the ideas by De Giorgi, Carriero and Leaci in the scalar case, the key point in obtaining 
Theorem~\ref{theop2} consists in establishing a one-sided Alfhors regularity for the jump set of (local) 
minimizers of the relaxed functional (\ref{e:weakMS}), also known in literature as \emph{density lower bound} 
estimate (for the precise formulation see Theorem~\ref{p:lowbound} below). 

In this perspective, we start off by considering the weak formulation of (\ref{eqgriffintro}). The functional
setting is provided by $SBD^2(\Omega)$, the space of fields $u\in L^1(\Omega,\R^n)$ with symmetrized 
distributional derivative $Eu:=(Du+Du^T)/2$ that is a bounded measure of the form
$Eu=e(u)\calL^n\res\Omega + [u]\odot\nu_u \calH^{n-1}\LL J_u$, with $[u]$ the jump of $u$, $J_u$ the 
$(n-1)$-rectifiable jump set of $u$, $\nu_u$ its normal, and with the properties $e(u)\in L^2(\Omega,\R^{n\times n})$ 
and $\calH^{n-1}(J_u)<\infty$. Here, $a\odot b=(a\otimes b+b\otimes a)/2$.
$SBD^2(\Omega)$ is a subset of the space of functions with bounded deformation $BD(\Omega)$. The latter was introduced and 
investigated in \cite{Suquet1978b,temam,TemamStrang1980,AnzellottiGiaquinta1980,KohnTemam1983} for the mathematical 
study of plasticity, damage and fracture models in a geometrically linear framework. 
Instead, $SBD^2(\Omega)$ provides the natural function space in the modeling of fracture in linear elasticity 
\cite{AmbrosioCosciaDalmaso1997,BellettiniCosciaDalmaso1998}.
Fine properties of $BD$ and $SBD^2$ are much less understood than those of their scalar counterparts $BV$ and 
$SBV^2$, respectively. Indeed, many standard technical tools are not available in this context, starting with basic ones
such as truncation results and the coarea formula. Despite this, recently several contributions have  improved
the understanding of such spaces 
\cite{chambolle,cha-gia-pon,gbd,ChambolleContiFrancfort,friedrich,friedrich1,ContiFocardiIurlano2016,DePhilippisRindler}. 

In view of the discussion above, the relaxation of (\ref{eqgriffintro}) is
\begin{equation}\label{eqgriffintroweak}
 E^*[u]:=\int_{\Omega} \Big(\frac12 \C e(u)\cdot e(u)+g(u)\Big) dx + \gamma\calH^{n-1}(J_u) 
\end{equation}
for $u\in SBD^2(\Omega)$. The density lower bound estimate for the jump set of
minimizers in this setting is the content of the ensuing theorem.
\begin{theorem}[Density lower bound]\label{p:lowbound}
If $u\in SBD^2(\Omega)\cap L^\infty(\Omega,\R^2)$ is a {minimizer} of the functional in (\ref{eqgriffintroweak})
with {$\|u\|_{L^\infty(\Omega,\R^2)}\leq M$}, then there exist $\vartheta_0$ and $R_0$, 
depending only on $\mathbb{C}$, {$g$}, $M$ and $\gamma$ such that if 
$0<\rho<R_0$, $x_0\in\overline{J_u}$, and $B_\rho(x_0)\subset\subset\Omega$, then 
\be{\label{e:lb}
\int_{B_\rho(x_0)}\frac12 \mathbb{C} e(u)\cdot e(u)dx+\gamma\calH^1(J_u\cap B_\rho(x_0))\geq \vartheta_0\rho.
}
Therefore,
\be{\label{e:saltochiuso}
\calH^1(\overline{J_u}\setminus J_u)=0.
}
\end{theorem}
Using this result, classical elliptic regularity yields that the minimizers $u$ belong to
$C^\infty(\Omega\setminus\overline{J_u},\R^2)$ (see for instance \cite{GiaquintaMartinazzi2012}), 
so that $(u,\overline{J_u})$ is a minimizer of the strong formulation of the problem in (\ref{eqgriffintro}).
This leads directly to the proof of Theorem \ref{theop2}.

The density lower bound estimate is a mild regularity result for the jump set of a minimizer $u$, therefore it is 
natural to analyze the infinitesimal behaviour of $u$ in points $x_0$ and, having selected a sequence 
$\rho_h\downarrow 0$, investigate the asymptotic of
\[
u_h(x):=\rho_h^{-\sfrac12}u(x_0+\rho_hx).
\]
We notice that the prefactor $\rho_h^{-\sfrac12}$ is needed to balance the different scaling of the volume and surface term in the 
energy ${E^*}$. Indeed, we have
\[
{E^*}[u_h;B_1]=\rho_h^{-1}{E^*}[u;B_{\rho_h}(x_0)],
\]
in the formula above the domains of integration are indicated explicitly.

The original proof of the density lower bound in formula (\ref{e:lb}) in the scalar case is indirect
\cite{DegiorgiCarrieroLeaci1989,AmbrosioFuscoPallara}. One
first constructs truncations of the rescaled functions $u_h$
and estimates them in $SBV$ using 
a Poincar\'e-Wirtinger type inequality. Then, using
Ambrosio's $SBV$ compactness theorem, one obtains convergence of a subsequence, 
and shows that the limit is a local minimizer of the bulk term of the Mumford-Shah energy $E^*_\mathrm{MS}$ 
restricted to Sobolev spaces, i.e., it is an harmonic function, and in particular smooth.
By a contradiction argument one then shows
that if in a ball the length of the jump set of a {minimizer $u$} of $E^*_\mathrm{MS}$ is sufficiently small and if {$u$ is not too far from being a (local) minimizer of the reduced functional}
\[u\mapsto\int_{\Omega} \frac12 |\nabla u|^2 dx + \gamma \calH^{n-1}(J_u) \,,\]
then in corresponding dyadic balls such an energy decays as fast as the Dirichlet integral for harmonic functions. 
From this, one deduces that the base point $x_0$ of the blow up process is not a jump point and the density lower bound follows at once.

The Poincar\'e-Wirtinger type inequality proven by De Giorgi, Carriero and Leaci states (in 2d) that
if $u\in SBV^2(B_1)$ and $\calH^{1}(J_u)$ is small, then there are $m\in \R$ and 
a modified function $\tilde u\in SBV^2(B_1)$ such that
$\|\tilde u-m\|_{L^{2}(B_1)}\le c\|\nabla u\|_{L^2(B_1,\R^2)}$, 
$|D\tilde u|(B_1)\le 2 \|\nabla u\|_{L^1(B_1,\R^2)}$, 
and $\tilde u=u$ on most of $B_1$, see \cite[Th. 4.14]{AmbrosioFuscoPallara} for details.
The function $\tilde u$ is obtained from $u$ by truncation, setting
$\tilde u(x)=\max\{\tau^-, \min\{u(x),\tau^+\}\}$ for suitable $\tau^\pm\in\R$ chosen through the coarea formula,
so that $\tilde u$ automatically fulfills $|D\tilde u|\le |Du|$. 
This procedure is not applicable to the vector-valued $BD$ functions which appear in the Griffith model,
since this space is not stable under truncation, and the coarea formula does not apply.

The key result to bypass such a problem is an approximation result for $SBD^p$ functions, $p\in(1,\infty)$, with small jump 
set with $W^{1,p}$ functions, stated below in the case of interest $p=2$ and established in \cite{ContiFocardiIurlanoRepr}.
This property yields an equivalent Poincar\'e-Wirtinger type inequality for $SBD^p$ functions, 
however restricted to two spatial dimensions. 

\begin{proposition}[Approximation of $SBD^2$ fields]\label{prop:ricopr}
There exist universal constants $c,\eta>0$ such that if $u\in SBD^2(B_{\rho})$, $\rho>0$, satisfies
\[
\Huno(J_u\cap B_{\rho})<{\eta\,(1-s)\frac\rho2}
\]
for some $s\in(0,1)$, then there are a countable family $\calF=\{B\}$ of closed balls of radius $r_B<(1-s)\rho/2$ 
with finite overlap, $\cup_\calF B \subset \subset B_\rho$ and a field $w\in SBD^2(B_{\rho})$  such that
\begin{enumerate}
\item $\rho^{-1}\sum_{\calF}{\Ln}\big(B\big)+\sum_{\calF}\Huno\big(\partial B\big)
\leq{\sfrac {c}\eta}\,\Huno(J_u\cap B_{\rho})$;
\item $\Huno\big(J_u\cap\cup_{\calF}\partial B\big)=\Huno\big((J_u\cap {B_{s\rho}})\setminus \cup_{\calF}B\big)=0$;
\item $w=u$ {$\Ln$-a.e.} on $B_{\rho}\setminus\cup_{\calF}B$;
\item {$w\in W^{1,2}(B_{s\rho},\Rdue)$} and $\Huno(J_w\setminus J_u)=0$;
\item If $u\in L^\infty(B_\rho,\R^2)$, then  $w\in L^\infty(B_\rho,\R^2)$ with $\|w\|_{L^\infty(B_\rho,\R^2)}\le \|u\|_{L^\infty(B_\rho,\R^2)}$;
\item 
\begin{equation}\label{e:volume}
 \int_{B}|e(w)|^2dx\leq c\int_{B}|e(u)|^2dx\qquad\text{for each $B\in\calF$};
\end{equation}
\item There is a skew-symmetric matrix $A$ such that
\begin{equation}\label{e:korn}
 \int_{B_{s\rho}\setminus\cup_{\calF}B}|\nabla u-A|^2dx\leq c\int_{B_{\rho}}|e(u)|^2dx.
\end{equation}
\end{enumerate}
\end{proposition}
Related results have been recently obtained in 
\cite{ChambolleContiFrancfort,friedrich,friedrich1}.

Proposition \ref{prop:ricopr} holds for any exponent $p$, is however restricted to two spatial dimensions. 
Its proof is based on covering the jump set of $u$ with balls such that the total length of the jump set contained in each of them 
is comparable (but significantly smaller than) the radius. Clearly, these balls cover a small part of $B_\rho$, and $u$ does not
need to be modified outside of them. In each of the balls, then, a new function is constructed by a finite-element approximation,
on a grid which refines close to the boundary of the ball, as was done in \cite{ContiSchweizer2006b} in the study of solid-solid 
phase transitions.
The key step is to show that one can choose such a grid with the 
property that all grid segments do not intersect the jump set of $u$. One then obtains an estimate of the oscillation of $u$
along the segments, and hence on the corners of each of the triangles which constitute the grid. Linear interpolation gives then
the desired extension. We refer to \cite{ContiFocardiIurlanoRepr} for the details of the proof.
In higher dimension, the same procedure would require finding a grid such that the edges do not intersect the jump set. This is however
not possible, at least with the strategy of \cite{ContiSchweizer2006b,ContiFocardiIurlanoRepr}, as was explained in those papers.

Proposition~\ref{prop:ricopr} is used in the proof of Theorem \ref{p:lowbound} to replace the truncation procedure
and the Poincar\'e-Wirtinger inequality. One then modifies the functions once more, subtracting not only a constant
as in the $BV$ case but also a linear function with skew-symmetric gradient, and obtains compactness. 
This permits to  classify the blow~up limits of
minimizers of (\ref{eqgriffintroweak}) with vanishing length of the jump set and to show that they minimize
a quadratic energy on Sobolev spaces, and therefore to conclude the proof of Theorem \ref{p:lowbound}.

\section*{Acknowledgements}
This work was partially supported 
by the Deutsche Forschungsgemeinschaft through the Sonderforschungsbereich 1060 
{\sl ``The mathematics of emergent effects''}.


\begin{thebibliography}{00}

\bibitem{AmbrosioCosciaDalmaso1997} 
L. Ambrosio, A. Coscia, G. Dal Maso, 
Fine properties of functions with bounded deformation, 
Arch. Rational Mech. Anal. 139 (1997), no. 3, 201--238.

\bibitem{amb-tort1}

L. Ambrosio, V.M. Tortorelli,
Approximation of functionals depending on jumps by elliptic functionals via $\Gamma$-convergence,
Comm. Pure Appl. Math. 43 (1990), n. 8, 999--1036.

\bibitem{AmbrosioFuscoPallara} 
L. Ambrosio, N. Fusco and D. Pallara,
Functions of Bounded Variation and Free Discontinuity Problems 
(Oxford University Press, Oxford, 2000).

\bibitem{AnzellottiGiaquinta1980}
G. Anzellotti, M. Giaquinta, 
Existence of the displacement field for an elastoplastic body subject to Hencky's law and von Mises 
yield condition, Manuscripta Math. 32 (1980), 1-2, 101--136.


\bibitem{BellettiniCosciaDalmaso1998} 
G. Bellettini, A. Coscia, G. Dal Maso, 
Compactness and lower semicontinuity properties in $SBD(\Omega)$, 
Math. Z. {228} (1998), no. 2, 337--351.

\bibitem{BourdinFrancfortMarigo2008} 
B. Bourdin, G. Francfort, J.-J. Marigo, 
The variational approach to fracture, 
Springer, New York, 2008. 

\bibitem{chambolle}

A. Chambolle, An approximation result for special functions with bounded
deformation,
J. Math. Pures Appl. (9) 83 (2004), n. 7, 929--954

\bibitem{ChambolleContiFrancfort}
A. Chambolle, S. Conti, G. Francfort,
Korn-Poincar\'e inequalities for functions with a small jump set,
Indiana Univ. Math. J., to appear; 
Preprint hal-01091710v1, 2014.

\bibitem{cha-gia-pon}
A. Chambolle, A. Giacomini, M. Ponsiglione,
Piecewise rigidity, J. Funct. Anal. 244 (2007), n. 1, 134--153.


\bibitem{ContiFocardiIurlano2016}
S. Conti and M. Focardi and F. Iurlano, 
Which special functions of bounded deformation have bounded variation?,
Preprint arXiv:1502.07464, 2015.

\bibitem{ContiFocardiIurlanoRepr} 
S. Conti and M. Focardi and F. Iurlano,
Integral representation for functionals defined on $SBD^p$ in dimension two,
Preprint arXiv:1510.00145, 2015.

\bibitem{ContiFocardiIurlano2016b}
S. Conti and M. Focardi and F. Iurlano,
Existence result for the static Griffith fracture model in two dimensions, in preparation.

\bibitem{ContiSchweizer2006b}
S. Conti and B. Schweizer,
Rigidity and $\Gamma$ convergence for solid-solid phase transitions with $SO(2)$-invariance,
Comm. Pure Appl. Math. 59,
(2006), 830--868.

\bibitem{gbd} G. Dal Maso, Generalised functions of bounded deformation,
J. Eur. Math. Soc. (JEMS) 15
(2013), n. 5,
1943--1997.

\bibitem{David2005}
G. David, Singular sets of minimizers for the Mumford-Shah
functional,
Progress in Mathematics
233,
Birkh\"auser Verlag, Basel, 2005.

\bibitem{DegiorgiCarrieroLeaci1989} 
E. De Giorgi, M. Carriero and A. Leaci,
Existence theorem for a minimum problem with free discontinuity set,  
Arch. Rational Mech. Anal. 108 (1989), no. 3, 195--218.

\bibitem{DePhilippisRindler}
G. De Philippis, F. Rindler, 
On the structure of $\mathcal{A}$-free measures and applications, 
Preprint arxiv-1601.06543, 2016.

\bibitem{FrancfortMarigo1998} 
G.A. Francfort and J.J. Marigo,
Revisiting brittle fractures as an energy minimization problem,
J. Mech. Phys. Solids 46 (1998), 1319--1342.

\bibitem{friedrich}
M. Friedrich,
A Korn type inequality in $SBD$ for functions with small jump sets,
Preprint arXiv:1505.00565, 2015.

\bibitem{friedrich1}
M. Friedrich,  
A Korn-Poincar\'e-type inequality for special functions of bounded deformation,
Preprint arxiv:1503.06755, 2015.

\bibitem{GiaquintaMartinazzi2012}
M. Giaquinta, L. Martinazzi,
An introduction to the regularity theory for elliptic systems, harmonic maps and minimal graphs,
Appunti Scuola Normale Superiore di Pisa (Nuova Serie),
11,
Edizioni della Normale, Pisa, 2012.

\bibitem{KohnTemam1983}
R. Kohn, R. Temam, 
Dual spaces of stresses and strains, with applications to Hencky plasticity,
Appl. Math. Optim. 10 (1983), n. 1, 1--35.

\bibitem{mum-shah} 
D. Mumford and J. Shah, 
Optimal approximation by piecewise smooth functions and associated variational problems, 
Comm. Pure Appl. Math. 17 (1989), 577--685.

\bibitem{Suquet1978b}
P.-M. Suquet, 
Existence et r\'egularit\'e des solutions des \'equations de la plasticit\'e,
C. R. Acad. Sci. Paris S\'er. A-B 286 (1978) A1201--A1204.

\bibitem{temam} 
R. Temam, 
Probl\`emes math\'ematiques en plasticit\'e, M\'ethodes 
Math\'ematiques de l'Informatique [Mathematical Methods of Information Science] 
12, Gauthier-Villars, Montrouge, 1983.

\bibitem{TemamStrang1980} 
R. Temam, G. Strang, 
Functions of bounded deformation, 
Arch. Rational Mech. Anal. 75 (1980/81), 7--21.

\end{thebibliography}
\end{document}